\documentclass{amsart}
\newcommand*\patchAmsMathEnvironmentForLineno[1]{%
  \expandafter\let\csname old#1\expandafter\endcsname\csname #1\endcsname
  \expandafter\let\csname oldend#1\expandafter\endcsname\csname end#1\endcsname
  \renewenvironment{#1}%
     {\linenomath\csname old#1\endcsname}%
     {\csname oldend#1\endcsname\endlinenomath}}%
\newcommand*\patchBothAmsMathEnvironmentsForLineno[1]{%
  \patchAmsMathEnvironmentForLineno{#1}%
  \patchAmsMathEnvironmentForLineno{#1*}}%
\AtBeginDocument{%
\patchBothAmsMathEnvironmentsForLineno{equation}%
\patchBothAmsMathEnvironmentsForLineno{align}%
\patchBothAmsMathEnvironmentsForLineno{flalign}%
\patchBothAmsMathEnvironmentsForLineno{alignat}%
\patchBothAmsMathEnvironmentsForLineno{gather}%
\patchBothAmsMathEnvironmentsForLineno{multline}%
}

\usepackage{amsfonts, amssymb, amsmath, latexsym, enumerate, array,caption,subcaption}
\usepackage{amscd}     
\usepackage[all]{xy}
\usepackage[pdftex]{graphicx}

\usepackage[latin1]{inputenc}
\usepackage{hyperref,cite,mdwlist,mathrsfs,color}

\usepackage[left=3.3cm,top=2.7cm,right=3.3cm]{geometry}

 \newtheorem*{pro}{Proposition}

\newtheorem{rem*}{Remark}
\newtheorem{rems*}{Remarks}
\newtheorem{problem*}{Problem}
\newtheorem{conj*}{Conjecture}

\newcommand{\sO}{\mathscr{O}}

\DeclareMathOperator{\HH}{H} \DeclareMathOperator{\hh}{h}

 \newcommand{\C}{\mathbb C}

 \newcommand{\p}{\mathbb P}

\SelectTips{cm}{} \newdir{ >}{{}*!/-5pt/@{>}}
\numberwithin{equation}{section}

\title[Eight cubes of linear forms in $\p^6$]{Eight cubes of linear forms in $\p^6$}

\author[G. Ilardi and J. Vall\`es]{Giovanna Ilardi and Jean Vall\`es}

\thanks{}

\begin{document}

\begin{abstract}
 Here we explain geometrically why the ideal $I=(L_1^3, \ldots, L_8^3)\subset \C[x_0,\ldots, x_6]$ has the WLP in degree $3$ and why it fails to have it in degree $5$.
\end{abstract}
\maketitle
\section{Introduction}
 In a private communication Rosa Mir\'o-Roig and Hoa Tran Quang informed us that 
 we wrongly affirm in \cite[Proposition 5.5]{DIV} that the ideal $I=(L_1^3, \ldots, L_8^3)$ fails the WLP in degree $3$, as it was conjectured in \cite[conjecture 6.6]{MMN1}.
 Indeed, computing explicitely the Hilbert functions of 
 the Artinian ring $A=\C[x_0,\ldots, x_6]/(L_1^3, \ldots, L_8^3)$ and $A/(L)$  
 where the $L_i$'s and $L$ are general linear forms, they observed that $I$ has the WLP in any degree except in degree $5$ meaning 
 that the multiplication map $\times L : A_i \rightarrow A_{i+1}$ has not maximal rank only when $i=5$. 
 
 \smallskip
 
 In this short note, we justify geometrically why $I$ has the WLP in degree $3$ and why it fails to have it in degree $5$.
 For the first one, our argument is based on a famous result by Alexander and Hirshowitz \cite{AH} who give a list of general sets of double points
 in $\p^n$ that do not impose independent conditions on hypersurfaces 
 of fixed degree. For the second one, we show that the failure is due to the existence of a pencil of cubics in $\p^5$ passing through $9$ quadruple  points in general position 
 in $\p^5$.
 
 \smallskip
 
 According to Mir\'o-Roig and Hoa Tran Quang degree $5$ is the only degree where there is a failure of the WLP. 
 Actually, if the failure can always be explained by a special geometric situation, 
 having the WLP, since it is expected in  general, is harder to prove. Indeed it would be necessary to have a theorem generalizing Alexander-Hirshowitz classification, that is
 a list of 
 general set of multiple points, with multiplicity bigger than $2$, that do not impose independent conditions on hypersurfaces 
 which is far to be known. 
 That's mainly why we do not propose a description in any degree.

\section{WLP in degree three}
\label{deg3}
Associated to $I$ there is the so-called Syzygy sheaf $K$ defined by:
$$
  \begin{CD}
 0@>>> K @>>>  \sO_{\p^{6}}^{8} @>(L_1^3,\ldots, L_8^3)>>
 \sO_{\p^{6}}(3) @>>> 0.
\end{CD}
$$
Let us recall that according to \cite{BK} we have $A_{d+i}=\HH^1(K(i))$.
By \cite[Lemma 5.2]{DIV} there are no syzygies of degree $\le 1$, i.e.  $\HH^0(K)=\HH^0(K(1))=0$.
Then taking the cohomology of the exact sequence defining $K$ we obtain the dimension of $A_3$ and $A_4$.
It occurs that the multiplication map 
$\times L: A_3\rightarrow A_4$ have maximal rank if and only if its cokernel has dimension exactly $78$. 
According to \cite{DIV} and in particular to \cite[Theorem 5.1]{DIV}, this cokernel coincide with the vector space of quartic cones in $\p^6$ with vertex at 
$\{L^{\vee}\}$ and $8$ double points at the $\{L_i^{\vee}\}$'s. But having a double point for a cone means that the line joining the double point to the vertex is also double. Then 
the dimension of this vector space of cones is exactly the dimension of the space of quartics in $\p^5$ with $8$ double points which is expected to be $126-6\times 8=126-48=78$. 
Let us point out that this expected dimension is actually the dimension since  $8$ double points in general position in $\p^5$ impose independent 
conditions to the quartics  because this is not one of the exceptional cases 
listed by Alexander and Hirshowitz in \cite{AH}.

\smallskip

This proves that
\begin{pro}
The ideal
$I=(L_1^3, \ldots, L_8^3)$ has the WLP in degree $3$ where
 $L_1, \ldots, L_8$ are general linear forms on $\p^6$.
\end{pro}
\begin{rem*}
 Since the multiplication map $\times L: A_3 \rightarrow A_4$ is injective we know, using \cite[Proposition 2.1]{MMN}, that 
 $\times L: A_i \rightarrow A_{i+1}$ are also injective for $i=0,1,2$.
\end{rem*}
\begin{rem*}
 The multiplication map $\times L: A_4 \rightarrow A_5$ is also injective according to Rosa Mir\'o-Roig and Hoa Tran Quang. If we want to apply the same technic than before we 
 have to compute the dimension of the space of quintics in $\p^5$ with $8$ triple points in general position. But we don't know if $8$ triple points in general position
 impose independent conditions on quintics and cannot conclude that the expected dimension is the true dimension.
\end{rem*}

\section{Failure of WLP in degree five}

We give now a geometric argument explaining why the map $\times L: A_5\rightarrow A_6$ is not injective.
Let us begin by computing the dimension of $A_5$ and $A_6$. 

Shifting by $2$ the exact sequence defining $K$ and computing the cohomology, we observe first that,  according to \cite[Lemma 5.2]{DIV},  $\HH^0(K(2))=0$.
This implies that $\mathrm{dim}(A_5)=\hh^1(K(2))=238$. Let us compute now 
$\mathrm{dim}(A_6)$. Shifting by $3$ the exact sequence defining $K$ and computing the cohomology, one finds: 
$$  \begin{CD}
 0@>>> \HH^0(K(3)) @>>>  \C^{672} @>>>
 \C^{924} @>>> A_6 @>>>0.
\end{CD}$$
The vector space $\HH^0(K(3))$ of syzygies of degree $3$ is not empty since it contains the Koszul relations, say $L_j^3. L_i^3+(-L_i^3).L_j^3=0$. These relations
are independent which gives $s:=\hh^0(K(3))\ge \binom{8}{2}=28$.
Then  $\mathrm{dim}A_6=924-672+s=280+(s-28)$\footnote{Actually we think that $s=28$ but we don't need the equality to prove the failure of WLP}.

\smallskip

Consequently the map $\times L:A_5 \rightarrow A_6$ is injective if and only if the cokernel has dimension $42+(s-28)$. According to \cite[Theorem 5.1]{DIV} and 
\cite[Theorem 13]{DI} the dimension of this cokernel is
the  sum of the dimension of the space of syzygies, that is $s$, and the dimension of the space of sextics cones in $\p^6$ with a vertex at $\{L^{\vee}\}$
and with $8$ quadruple points
at the $\{L_i^{\vee}\}$'s. As we wrote before in section (\ref{deg3}) these sextics correspond to sextics in $\p^5$ with $8$ quadruple points in general position. 
Its expected dimension is $14$ which added to the $s$ syzygies would give $s+14=42+(s-28)$.

\smallskip

But a special situation occurs: through $9$ double points in general position in $\p^5$ there is a pencil of cubics, let say $(C_1,C_2)$. There are also $8$ cubics through
our $8$ points among these $9$, let say $(C_1,C_2,\ldots, C_8)$ since of course $C_1$ and $C_2$ belong to this linear system. Now the vector space 
$$ (C_1^2, C_1C_2, \ldots, C_1C_8, C_2^2, C_2C_3, \ldots, C_2C_8)$$ has dimension $15$ and it consists in sextics with $8$ quadruple  points. This proves that the cokernel 
of $A_5 \rightarrow A_6$ has dimension at least $s+15$. This proves that the ideal $I$ fails the WLP in degree $5$.

\begin{rem*}
 We observe  that $9$ quadruple points in general position does not impose independent conditions on sextics of $\p^5$. Indeed if we compute directly there is no such sextics.
 But of course the existence of a pencil of cubics passing through these $9$ double points give sextics with $9$ quadruple points
 which proves that the dimension is strictly bigger than the expected one.
\end{rem*}

\end{document}